\documentclass[12pt]{article}

\usepackage{amsmath, amsthm, amsfonts, amssymb}

\newtheorem{thm}{Theorem}[section]

\newtheorem{lem}[thm]{Lemma}

\newtheorem{rem}[thm]{Remarque}

\newcommand{\dem}{{\bf Proof.\ }}
\newcommand{\fin}{\rule{0.5em}{0.5em}}
\newcommand{\semi}{$\left\{T(t)\right\}_{t\geq 0}\;$}
\newcommand{\semia}{$\left\{T^{*}(t)\right\}_{t\geq 0}\;$}

\title{{\bf Uniqueness of a pre-generator for $C_0$-semigroup on a general locally convex vector space}\thanks{This
work is partially supported by {\it Yangtze Research Programme},
Wuhan University, China, and the {\it Town Council} of Hunedoara,
Romania.}}
\author
{Ludovic Dan LEMLE\thanks{Laboratoire de Math\'ematiques, CNRS-UMR 6620, Universit\'e
Blaise Pascal, 63177 Aubi\`ere, France, e-mail: {\tt lemle.dan@fih.upt.ro}\quad and Engineering Faculty, "Politehnica"
University, 331128 Hunedoara, Romania.}\quad\quad Liming
WU\thanks{Laboratoire de Math\'ematiques, CNRS-UMR 6620,
Universit\'e Blaise Pascal, 63177 Aubi\`ere, France. e-mail: {\tt
Li-Ming.Wu@math.univ-bpclermont.fr}\quad and Department of
Mathematics, Wuhan University, 430072-Hubei, P.R. China.}}

\date{version 25 January 2007}

\def\D{\cal{D}}

\def\N{\mathbb{N}}

\def\R{\mathbb{R}}

\def\X{\cal{X}}
\def\Y{\cal{Y}}

\begin{document}

\maketitle

\begin{abstract}
\noindent The main purpose is to generalize a theorem of Arendt about uniqueness of $C_0$-semigroups from Banach space setting to the general locally convex vector spaces, more precisely, we show that cores are the only domains of uniqueness for
$C_0$-semigroups on locally convex spaces. As an application, we find a necessary and sufficient condition for that the mass transport equation has one unique $L^1(\R^d,dx)$ weak solution.
\end{abstract}

{\bf Running head:} Uniqueness of a pre-generator

\section{Framework and main result}
The theory of $C_0$-semigroups of linear operators in Banach spaces
was extended by {\sc Schwartz} \cite{schwartz'58}, {\sc Miyadera}
\cite{miyadera'59}, {\sc Yosida} \cite{yosida'71}, {\sc Komatsu}
\cite{komatsu'64}, {\sc K\~omura} \cite{komura'68} and others to the
case of equicontinuous $C_0$-semigroups of linear operators in
locally
convex spaces. Let $(\X,\beta)$ be a locally convex Hausdorff space. Recall that a
family \semi of linear continuous operators on $\X$ is called a {\it $C_0$-semigroup on $\X$} if the following
properties holds:\\
(i) $T(0)=I$;\\
(ii) $T(t+s)=T(t)T(s)$, for all $t,s\geq 0$;\\
(iii) $\lim_{t\searrow 0}T(t)x=x$, for all $x\in\X$;\\
(iv) there exist a number $\omega\in\R$ such that the family
$\left\{e^{-\omega t}T(t)\right\}_{t\geq 0}$ is equicontinuous.\\
Furthermore we say that \semi is an equicontinuous $C_0$-semigroup if $\omega=0$ in (iv). 
The equicontinuity must be considered in the sense of seminorms: a
family $\cal F$ of linear operators on $\X$ is said to be
equicontinuous if for each continuous seminorm $p$ on $\X$, there is
a continuous seminorm $q$ on $\X$ such that
$$
p(Tx)\leq q(x)\quad,\quad\forall T\in{\cal F}\mbox{ and }\forall
x\in \X.
$$
The {\it infinitesimal generator} of $C_0$-semigroup \semi is a
linear operator $\cal L$ defined on the domain
$$
{\cal D}({\cal L})=\left\{x\in{\X}\:\left|\:\lim_{t\searrow 0}\frac{T(t)x-x}{t}\mbox{ exists in }(\X,\beta)\right.\right\}
$$
by
$$
{\cal L}x=\lim_{t\searrow 0}\frac{T(t)x-x}{t}\quad,\quad\forall
x\in{\cal D}({\cal L}).
$$
If the locally convex Hausdorff space $(\X,\beta)$ is assumed to be
sequentially complete, then:\\
(i) $\cal L$ is a densely defined and closed operator;\\
(ii) the resolvent $R(\lambda;{\cal L})=(\lambda I-{\cal L})^{-1}$,
for any $\lambda\in\rho({\cal L})$ (the resolvent set of $\cal L$)
is well defined, continuous on $\X$ and satisfies the equality
$$
R(\lambda;{\cal L})=\int\limits_{0}^{\infty}\!e^{-\lambda
t}T(t)x\:dt\quad,\quad\forall \lambda>\omega\mbox{ and }\forall
x\in\X.
$$
Let ${\cal A}:{\X}\longrightarrow{\X}$ be a linear operator with
domain $\D$ dense in $\X$. $\cal A$ is said to be a {\it
pre-generator}, if there exists some $C_0$-semigroup on $\X$ such
that its generator $\cal L$ extends $\cal A$. We say that $\cal A$
is an {\it essential generator} in $\X$ (or $\X$-unique), if $\cal
A$ is closable and its closure $\overline{\cal A}$ with respect to
$\beta$ is the generator of some $C_0$-semigroup on $\X$.\\
In general, for a $C_0$-semigroup \semi on $({\X},\beta)$, its
adjoint semigroup \semia is  no longer strongly continuous on the
dual topological space ${\Y}$ of $(\X,\beta)$ with respect to the
strong topology $\beta({\Y},{\X})$ of $\Y$. In
\cite[p.563]{wu-zhang'06} {\sc  Zhang} and the second named author introduced a new
topology on $\Y$ for which the usual semigroups in the literature
becomes $C_0$-semigroups. That is {\it the topology of uniform
convergence on compact subsets of $({\X},\beta)$}, denoted by ${\cal
C}({\Y},{\X})$. If moreover, $({\X},\beta)$ is assumed to be
quasi-complete (i.e., the bounded and closed subsets of $(\X,\beta)$
are complete) then $(\Y,{\cal C}(\Y,\X))^{*}=(\X,\beta)$ and if
\semi is a $C_0$-semigroup on $({\X},\beta)$ with generator $\cal
L$, then \semia is
a $C_0$-semigroup on $({\Y},{\cal C}({\Y},{\X}))$ with generator ${\cal L}^{*}$.\\
The main purpose of this paper is to furnish a proof for the
difficult implication of a theorem of {\sc Wu} and {\sc Zhang}
\cite[Theorem 2.1, p. 570]{wu-zhang'06} concerning uniqueness of
pre-generators on locally convex spaces (first time formulated in \cite[Remarks (2.v), p. 292]{wu'98}).
\begin{thm}\label{1.1}
Let $(\X,\beta)$ be a locally convex Hausdorff sequentially complete
space and $\cal A$ a linear operator on $\X$ with domain $\D$ (the
test-function space) which is assumed to be dense in $(\X,\beta)$.
Assume that there is a $C_0$-semigroup \semi on $(\X,\beta)$ such
that its generator $\cal L$ is an extension
of $\cal A$ (i.e., $\cal A$ is a pre-generator in $(\X,\beta)$). The following assertions are equivalents:\\
(i) $\cal A$ is a $\X$-essential generator (or $\X$-unique);\\
(ii) the closure of $\cal A$ in $(\X,\beta)$ is exactly $\cal L$ (i.e., $\D$ is a core for $\cal L$);\\
(iii) ${\cal A}^{*}={\cal L}^{*}$ which is the generator of the dual $C_0$-semigroup \semia on
$({\Y},{\cal C}({\Y},{\X}))$;\\
(iv) for some $\lambda>\omega$ ($\omega\in\R$ is the constant in
definition of $C_0$-semigroup \semi),
the range $(\lambda I-{\cal A})(\D)$ is dense in $(\X,\beta)$;\\
(v) (Liouville property) for some $\lambda>\omega$, ${\cal
K}\mbox{\it er }(\lambda I-{\cal A}^{*})=\{0\}$
 (i.e., if $y\in{\cal D}({\cal A}^{*})$ satisfies $(\lambda I-{\cal A}^{*})y=0$, then $y=0$);\\
(vi) (uniqueness of solutions for the resolvent equation) for all
$\lambda>\omega$ and all $y\in{\Y}$, the resolvent equation of
${\cal A}^{*}$
$$
(\lambda I-{\cal A}^{*})z=y
$$
has the unique solution $z=((\lambda I-{\cal L})^{-1})^{*}y=(\lambda I-{\cal L}^{*})^{-1}y$;\\
(vii) (uniqueness of strong solutions for the Cauchy problem) for
each $x\in{\cal D}\left(\overline{\cal A}\right)$, there is a
$(\X,\beta)$-unique strong solution $v(t)=T(t)x$ of the Cauchy
problem (or the Kolmogorov backward equation)
$$
\left\lbrace\begin{array}{l}
\partial_tv(t)=\overline{\cal A}v(t)\\
v(0)=x
\end{array}
\right.
$$
i.e., $t\mapsto v(t)$ is differentiable from $\R^{+}$ to
$(\X,\beta)$ and its derivative $\partial_tv(t)$ coincides
with $\overline{\cal A}v(t)$;\\
(viii) (uniqueness of weak solutions for the dual Cauchy problem)
for every $y\in{\Y}$, the dual Cauchy problem (or the Kolmogorov
forward equation)
$$
\left\lbrace\begin{array}{l}
\partial_tu(t)={\cal A}^{*}u(t)\\
u(0)=y
\end{array}
\right.
$$
has a $(\Y,{\cal C}(\Y,\X))$-unique weak solution $u(t)=T^{*}(t)y$.
More precisely, there is a unique function $\R^{+}\ni t\longmapsto
u(t)=T^{*}(t)y$ which is continuous from $\R^{+}$ to $({\Y},{\cal
C}({\Y},{\X}))$ such that
$$
\left\langle x,u(t)-y\right\rangle=\int\limits_{0}^{t}\!\left\langle
{\cal A}x,u(s)\right\rangle ds\quad,\quad\forall x\in{\D};
$$
(ix) there is only one $C_0$-semigroup on $\X$ such that its
generator extends $\cal A$.
\end{thm}
Many equivalence relations above, especially the equivalence between (i), (vii), (viii) and (ix), are fundemental and well studied in the Banach space setting, see {\sc Arendt} \cite{arendt'86}, {\sc Eberle} \cite{eberle'97}, {\sc Djellout} \cite{djellout'97} and the second named author \cite{wu'98} and \cite{wu'99}, etc. In the local convex space framework, the equivalences between (i)-(viii) are proved in {\sc Zhang} and the second named author \cite{wu-zhang'06}.\\
The main purpose of this paper is to prove the equivalence between (i) and (ix). This equivalence in the Banach space setting is the well known Arendt theorem. The implication (i)$\Rightarrow$(ix) is immediate. Indeed, if $\cal A$ is an
essential generator in $\X$ and if we suppose that $\cal L$ and
$\cal L'$ are generators of some $C_0$-semigroups which extends
$\cal A$, then by equivalence of (i)$\Leftrightarrow$(ii), we have
${\cal L}=\overline{\cal A}={\cal L'}$. It follows that there is
only one $C_0$-semigroup
on $\X$ such that its generator extends $\cal A$.\\
The sufficiency of (ix) is difficult. We shall follow the strategy of Arendt in the Banach space setting, but several basic ingredients require much more difficult proofs in the actual locally convex vector space setting. The main idea for overcoming those difficulties is to use the notion of calibration and to choose a "good" calibration.

This paper is organized as follows. In the next section we show that
cores are the only domains of uniqueness for $C_0$-semigroups on
locally convex spaces and we prove the sufficiency of (ix) in
Theorem \ref{1.1}. Finally is presented the $L^1$-uniqueness for weak solution of mass transport equation.

\section{Domains of uniqueness}
Recall at first several well known facts for calibration. A {\it
calibration} for a locally convex space $(\X,\beta)$ is a family
$\Gamma$ of continuous seminorms which induces the topology $\beta$
of $\X$. Such a family of seminorms was used by {\sc Fattorini}
\cite{fattorini'68}, {\sc Moore} \cite{moore'69}, {\sc
Chilana} \cite{chilana'70}, {\sc Choe} \cite{choe'85} and others.\\
Let $p\in\Gamma$. A linear operator $T$ on $\X$ is said to be {\it
$p$-continuous} if
$$
\tilde{p}(T):=\sup_{p(x)\leq 1}p(Tx)<\infty
$$
and is said to be {\it $\Gamma$-continuous} if it is $p$-continuous
for every $p\in\Gamma$. We say that $T$ is {\it $\Gamma$-bounded} if
$$
\|T\|_{\Gamma}:=\sup_{p\in\Gamma}\tilde{p}(T)<\infty\quad.
$$
If $\|T\|_{\Gamma}\leq 1$, then we say that $T$ is a {\it
$\Gamma$-contraction}.\\
The following result obtained by {\sc Moore} \cite[Theorem 4, p.
70]{moore'69}, give a very nice characterisation of equicontinuous
semigroups.
\begin{lem}\label{2.1}
A semigroup $\cal F$ of linear operators on $\X$ is equicontinuous
if and only if there is a calibration $\Gamma$ for $\X$ such that
$\cal F$ is a semigroup of $\Gamma$-contraction.
\end{lem}
Finally, the following perturbation result of {\sc Choe}
\cite[Corollary 5.4, p. 312]{choe'85}, plays a key role in the proof
of our next theorem.
\begin{lem}\label{2.2}
Let $\Gamma$ be a calibration for a locally convex space
$(\X,\beta)$. If $A$ is the generator of a $C_0$-semigroup on $\X$
and $B$ is a $\Gamma$-bounded linear operator on $\X$, then $A+B$ is
the generator of a $C_0$-semigroup on $\X$.
\end{lem}
We turn now to the job. We begin with the following theorem which is
well known in the Banach space setting (see \cite[Theorem 1.33, p.
46]{arendt'86}).
\begin{thm}\label{2.3}
Let $(\X,\beta)$ be a locally convex Hausdorff sequentially complete
space, $\left\{T(t)\right\}_{t\geq 0}$ a $C_0$-semigroup on $\X$
with generator $\cal L$ and $\D$ a subspace of ${\cal D}({\cal L})$.
Consider the restriction $\cal A$ of $\cal L$ to $\D$. If $\D$ is
not a core of $\cal L$, then there exists an infinite number of
extensions of $\cal A$ which are generators.
\end{thm}
\dem {\bf Step 1.} Endow ${\cal D}({\cal L})$ with the graph topology $\beta_{\cal
L}$ of ${\cal L}$ induced by the $\beta$-topology. If in contrary
$\D$ is not a core of $\cal L$, then $\D$ is not dense in ${\cal
D}({\cal L})$ with respect to the graph topology $\beta_{\cal L}$.
By Hahn-Banach theorem there exist some non-zero
linear functional $\phi$ continuous on ${\cal D}({\cal L})$ with
respect to the graph topology $\beta_{\cal L}$ such that $\phi(x)=0$
for all $x\in{\D}$. Fix some $u\in{\cal D}({\cal L})$, $u\neq 0$, we
consider the linear operator
$$
C:{\cal D}({\cal L})\longrightarrow{\cal D}({\cal L})
$$
$$
Cx=\phi(x)u\quad,\quad\forall x\in{\cal D}({\cal L}).
$$
Then $C$ is $\beta_{\cal L}$-continuous (i.e. continuous with respect to the graph topology $\beta_{\cal L}$)
on ${\cal D}({\cal L})$. Notice that $C$ is $\beta_{\cal L}$-continuous iff for some (or all)
$\lambda_0\in\rho({\cal L})$
$$
\tilde{C}:=(\lambda_0I-{\cal L})CR(\lambda_0;{\cal L})
$$
is $\beta$-continuous on $\X$.\\
Let $\Theta=CR(\lambda_0;{\cal L})$. Since for all $x\in\X$ we have
$$
\Theta x=CR(\lambda_0;{\cal L})x=\phi\left(R(\lambda_0;{\cal
L})x\right)u
$$
$$
\Theta^2x=\Theta(\Theta x)=\phi\left(R(\lambda_0;{\cal L})\Theta x\right)u=
\phi\left(R(\lambda_0;{\cal L})\phi\left(R(\lambda_0;{\cal
L})x\right)u\right)u=
$$
$$
=\phi\left(R(\lambda_0;{\cal L})x\right)\phi\left(R(\lambda_0;{\cal
L})u\right)u,
$$
and successively
$$
\Theta^nx=\phi\left(R(\lambda_0;{\cal L})x\right)\phi^{n-1}\left(R(\lambda_0;{\cal
L})u\right)u 
$$
for all $n\in\N^{*}$. One can  take $u\in{\cal D}({\cal L})$, $u\neq 0$
such that
$$
\left|\phi\left(R(\lambda_0;{\cal L})u\right)\right|<1.
$$
Therefore the linear operator $U=I-CR(\lambda_0;{\cal L})$ is
invertible and both $U$ and $U^{-1}$,
$$
U^{-1}x=\sum\limits_{n=0}^{\infty}\Theta^nx=x+\phi\left(R(\lambda_0;{\cal L})x\right)\frac{1}{1-\phi\left(R(\lambda_0;{\cal L})u\right)}u
$$
are $\beta$-continuous on $\X$.
Moreover, as in the proof of \cite[Theorem 1.31, p. 45]{arendt'86},
we have
\begin{eqnarray*}
U\left({\cal L}+\tilde{C}\right)U^{-1}&=&U\left({\cal L}-\lambda_0 I+\lambda_0 I+\tilde{C}\right)U^{-1}=\\
&=&U\left({\cal L}-\lambda_0 I+\tilde{C}\right)U^{-1}+\lambda_0 I=\\
&=&U\left({\cal L}-\lambda_0 I+\left(\lambda_0 I-{\cal L}\right)CR(\lambda_0;{\cal L})\right)U^{-1}+\lambda_0 I=\\
&=&U\left({\cal L}-\lambda_0 I\right)\left(I-CR(\lambda_0;{\cal L})\right)U^{-1}+\lambda_0 I=\\
&=&U\left({\cal L}-\lambda_0 I\right)+\lambda_0 I=\\
&=&\left[I-CR(\lambda_0;{\cal L})\right]\left({\cal L}-\lambda_0I\right)+\lambda_0 I=\\
&=&{\cal L}-\lambda_0 I+C+\lambda_0 I={\cal L}+C\quad.
\end{eqnarray*}
Now we have only to prove that ${\cal L}+\tilde{C}$ is the generator
of some $C_0$-semigroup on $(\X,\beta)$.\\
{\bf Step 2.} To apply Lemma \ref{2.2} of Choe, we have to find a good calibration, which is the main difficult point. Since
$\left\{T(t)\right\}_{t\geq 0}$ is a $C_0$-semigroup on $\X$, there
exists a number $\omega\in\R$ such that $\left\{e^{-\omega
t}T(t)\right\}_{t\geq 0}$ is equicontinuous. According to Lemma
\ref{2.1}, there is a calibration $\Gamma$ for $(\X,\beta)$ such
that
$$
\left\|e^{-\omega t}T(t)\right\|_{\Gamma}\leq 1\quad,\quad\forall t\geq 0.
$$
For each $p\in\Gamma$ we define
$$
\hat{p}(x)=\sup_{t\geq 0}\left[p\left(e^{-\omega
t}T(t)x\right)+\left|\phi\left(R(\lambda_0;{\cal
L})e^{-\omega
t}T(t)x\right)\right|\right]\quad,\quad\forall x\in{\X}.
$$
As $\hat{p}\geq p$ and $\hat{p}$ is continuous, the family
$\hat{\Gamma}=\left\{\hat{p}\:|\:p\in\Gamma\right\}$ is
another calibration of $(\X,\beta)$, which will be our calibration. We consider now the
$\hat{\Gamma}$-norm
$$
\left\|\tilde{C}\right\|_{\hat{\Gamma}}=\sup_{\hat{p}\in\hat{\Gamma}}\sup_{\hat{p}(x)\leq
1}\hat{p}\left(\tilde{C}x\right)
$$
and we prove that $\tilde{C}$ is $\hat{\Gamma}$-bounded, i.e.
$$
\left\|\tilde{C}\right\|_{\hat{\Gamma}}<\infty\quad.
$$
Let $x\in\X$. Then we have
\begin{eqnarray*}
\tilde{C}x&=&\left(\lambda_0 I-{\cal L}\right)CR(\lambda_0;{\cal L})x=\left(\lambda_0I-{\cal L}\right)\phi\left(R(\lambda_0;{\cal L})x\right)u=\\
&=&\phi\left(R(\lambda_0;{\cal L})x\right)\left(\lambda_0I-{\cal
L}\right)u=\phi\left(R(\lambda_0;{\cal L})x\right)v
\end{eqnarray*}
where we denote $\left(\lambda_0I-{\cal L}\right)u=v$. Therefore
\begin{eqnarray*}
\hat{p}\left(\tilde{C}x\right)&=&\hat{p}\left(\phi\left(R(\lambda_0;{\cal L})x\right)v\right)=\left|\phi\left(R(\lambda_0;{\cal L})x\right)\right|\hat{p}(v)\\
&\leq&\left[p(x)+\left|\phi\left(R(\lambda_0;{\cal L})x\right)\right|\right]\hat{p}(v)\leq\hat{p}(x)\hat{p}(v)\quad.
\end{eqnarray*}
Consequently
$$
\sup_{\hat{p}(x)\leq 1}\hat{p}\left(\tilde{C}x\right)\leq\sup_{\hat{p}(x)\leq
1}\hat{p}(x)\hat{p}(v)\leq\hat{p}(v)\quad.
$$
Then $\left\|\tilde{C}\right\|_{\hat{\Gamma}}<\hat{p}(v)$, i.e. $\tilde{C}$ is $\hat{\Gamma}$-bounded. So by Lemma \ref{2.2} ${\cal L}+\tilde{C}$ generate a $C_0$-semigroup $\left\{\tilde{S}(t)\right\}_{t\geq 0}$. Consequently
$$
S(t)=U\tilde{S}(t)U^{-1}
$$
is a $C_0$-semigroup whose generator is ${\cal L}+C$ and ${\cal L}+C/_{\cal D}={\cal L}/_{\cal D}$. As the choice of $u$ above is infinite, we have proved the result. \fin\\
{\bf Proof of Theorem \ref{1.1} (ix)$\Rightarrow$(i)} Suppose that
there is only one $C_0$-semigroup on $\X$ such that its generator
extends $\cal A$. By the Theorem \ref{2.3} it follows that ${\D}$ is
a core of ${\cal L}$. Therefore $\overline{{\cal L}/{_{\D}}}={\cal
L}$. But ${\cal A}={\cal L}/{_{\D}}$, we conclude that
$\cal A$ is a $\X$-unique. \fin
\begin{rem}
\em If $\cal A$ is a second order elliptic differential operator with ${\cal D}=C_0^\infty(D)$, then the weak solutions for
the dual Cauchy problem in the Theorem \ref{1.1} (viii) correspond
exactly to those in the distribution sense in the theory of partial
differential equations and the dual Cauchy problem becomes the Fokker-Planck equation for heat diffusion. We must remarks the important equivalences
between the $\X$-uniqueness of the linear operator $\cal A$, the
$\X$-uniqueness of strong solutions for the Cauchy problem and the
$\Y$-uniqueness of weak solutions for the dual Cauchy problem
associated with $\cal A$.
\end{rem}

\section{$L^1(\R^d,dx)$-uniqueness of weak solutions for the mass transport equation}

Consider the operator
$$
{\cal A}f=b\nabla f\quad,\quad\forall f\in
C_0^\infty\left(\R^d\right)
$$
where the vector field $b:\R^d\rightarrow\R^d$ is locally Lipschitz. Let $\partial$ be the point at infinity of $\R^d$. Consider the ordinary differential equation (ODE)
$$
\left\lbrace\begin{array}{l}
dX_t=b(X_t)dt\\
X(0)=x
\end{array}
\right.\quad.
$$
For every $x\in\R^d$, there is a unique solution
$(X_t(x))_{0\leq t<\tau_e}$, where
$$
\tau_e=\inf\left\{t\geq 0\:|\:X_t=\partial\right\}
$$
is the explosion time. Then the family $\left\{P_t\right\}_{t\geq 0}$, where
$$
P_tf(x)=f(X_t(x))1_{[t<\tau_e]}
$$
is a $C_0$-semigroup on $L^\infty\left(\R^d,dx\right)$ with respect
to the topology ${\cal C}\left(L^\infty,L^1\right)$ and
$$
f(X_t)-f(X_0)=\int\limits_0^t\!b\nabla f(X_s)ds\quad,\quad\forall f\in{\cal C}_0^\infty(\R^d).
$$
Therefore $f$ belongs to the domain of the generator ${\cal
L}_{(\infty)}$ of $C_0$-semigroup $\left\{P_t\right\}_{t\geq 0}$ on
$L^\infty\left(\R^d,dx\right)$ and 
$${\cal L}_{(\infty)}f={\cal A}f=b\nabla f\quad.
$$
Consequently, $\left({\cal
A},C_0^\infty(\R^d)\right)$ is a pre-generator on
$\left(L^\infty\left(\R^d,dx\right),{\cal
C}\left(L^\infty,L^1\right)\right)$. So we can study the
$\left(L^\infty\left(\R^d,dx\right),{\cal
C}\left(L^\infty,L^1\right)\right)$-uniqueness of the operator
$\left({\cal A},C_0^\infty(\R^d)\right)$.\\
Consider at first the one-dimensional operator
$$
{\cal A}f=bf^{'}\quad,\quad\forall f\in C_0^\infty(\R)
$$
where $b$ is a locally Lipschitz continuous function on $\R$ such
that $b(x)>0$, for all $x\in\R$.
Let ${\cal A}^{*}:{\cal D}\left({\cal A}^{*}\right)\subset
L^1(\R,dx)\rightarrow L^1(\R,dx)$ the adjoint operator of $\cal A$
and $h\in L^1(\R,dx)$ such that $h\in{\cal D}\left({\cal
A}^{*}\right)$ and
$$
{\cal A}^{*}h=\lambda h\quad.
$$
Then $h$ solve the ODE in the distribution sense
$$
-(bh)^{'}=\lambda h\quad.
$$
Then $bh$ is absolutely continuous where it follows that $h$ is absolutely continuous in the set $\{x\in\R\:|\:b(x)\neq 0\}$.
\begin{thm}
Assume that $b$ is locally Lipschitzian and $b(x)>0$ over $\R$. Then $\left({\cal A},C_0^\infty(\R)\right)$ is
$\left(L^\infty\left(\R,dx\right),{\cal
C}\left(L^\infty,L^1\right)\right)$-unique if and only if
$$
\int\limits_{-\infty}^0\!\frac{1}{b(x)}dx=\infty\quad.
$$
\end{thm}
\dem As shown before $\left({\cal A},C_0^\infty(\R)\right)$ is a
pre-generator on $\left(L^\infty\left(\R,dx\right),{\cal
C}\left(L^\infty,L^1\right)\right)$.\\
{\it Sufficiency.} Suppose in contrary that $\left({\cal A},C_0^\infty(\R)\right)$ is
not $\left(L^\infty\left(\R,dx\right),{\cal
C}\left(L^\infty,L^1\right)\right)$-unique. Then there is a function
$h\in L^1(\R,dx)$, $h\neq 0$ such that
$$
\left(I-{\cal A}^{*}\right)h=0
$$
in the sense of distributions. Thus we may assume that $h$ is itself absolutely continuous and solves the ODE
$$
-(bh)^{'}=h\quad.
$$
Then
$$
h^{'}=-\frac{1+b^{'}}{b}h\quad,
$$
we find
$$
h(x)=h(0)e^{-\int\limits_0^x\!\frac{1+b^{'}(s)}{b(s)}ds}=h(0)\frac{b(0)}{b(x)}e^{-\int\limits_0^x\frac{1}{b(s)}ds}.
$$
Because $h\neq 0$, we have $h(0)\neq 0$. Then
$$
\int\limits_{\R}|h(x)|dx=|h(0)||b(0)|\int\limits_{\R}\frac{1}{b(x)}e^{-\int\limits_0^x\frac{1}{b(s)}ds}dx
$$
and for
$$
u(x)=\int\limits_0^x\frac{1}{b(s)}ds
$$
$u(-\infty)=-\infty$, we obtain
$$
\int\limits_{\R}|h(x)|dx=|h(0)||b(0)|\int\limits_{u(-\infty)}^{u(\infty)}e^{-u}du=\infty
$$
which is in contradiction with the assumption that $h\in
L^1(\R,dx)$.\\ 
{\it Necessity.} If in contrary
$$
\int\limits_{-\infty}^0\!\frac{1}{b(x)}dx<\infty\quad,
$$
then
$$
h(x)=\frac{b(0)}{b(x)}e^{-\int\limits_{0}^x\!\frac{1}{b(s)}\:ds}\in L^1(\R,dx)
$$
and
$$
(I-{\cal A}^{*})h=0
$$
which is contradictory to the fact that $\cal A$ is $\left(L^\infty\left(\R,dx\right)\right)$-unique. \fin\\
We can formulate next
\begin{thm}
Let $b$ be a locally Lipschitz continuous function on $\R$ such that
for some $c_0<c_N\in\R$, $b1_{(-\infty,c_0]}$ and
$b1_{[c_N,+\infty)}$ keep a constant sign (non-zero). Then $\left({\cal A},C_0^\infty(\R)\right)$ is
$\left(L^\infty\left(\R,dx\right),{\cal
C}\left(L^\infty,L^1\right)\right)$-unique if and only if
$$
\int\limits_{-\infty}^{c_0}\frac{1}{b^{+}(x)}dx=\int\limits_{c_N}^{+\infty}\frac{1}{b^{-}(x)}dx=+\infty.
$$
\end{thm}
\dem
{\it Sufficiency.} Suppose in contrary that $\left({\cal A},C_0^\infty(\R)\right)$ is
not $\left(L^\infty\left(\R,dx\right),{\cal
C}\left(L^\infty,L^1\right)\right)$-unique. Then there is a function
$h\in L^1(\R,dx)$, $h\neq 0$ such that
$$
\left(I-{\cal A}^{*}\right)h=0
$$
in the sense of distributions. Then $bh$ is absolutely continuous over $\R$.\\
Because $b1_{(-\infty,c_0]}$ and
$b1_{[c_N,+\infty)}$ keep a constant sign, we may suppose that
$$
\left\{x\in\R\:|\:b(x)=0\right\}=\{x_1<x_2<...<x_N\}\subset[c_0,c_N].
$$
{\bf Step 1.} Let $x\in I_k=(x_k,x_{k+1})$ and $c_k\in I_k$, $k\in\{1,2,...,N-1\}$. Since $h$ is absolutely continuous over $I_k$, we have
$$
h(x)=h(c_k)\frac{b(c_k)}{b(x)}e^{-\int\limits_{c_k}^{x}\!\frac{1}{b(s)}\:ds}\quad,\quad\forall x\in I_k.
$$
Because $h\neq 0$ in $I_k$, we have
\begin{itemize}
\item if $b(x)>0$ for all $x\in I_k$, then 
$$
\lim_{x\searrow x_k}b(x)h(x)=b(c_k)h(c_k)e^{\int\limits_{x_k}^{c_k}\!\frac{1}{b(s)}\:ds}\mbox{ not exist}
$$
\item if $b(x)<0$ for all $x\in I_k$, then 
$$
\lim_{x\nearrow x_{k+1}}b(x)h(x)=b(c_k)h(c_k)e^{-\int\limits_{c_k}^{x_{k+1}}\!\frac{1}{b(s)}\:ds}\mbox{ not exist}
$$
\end{itemize}
But all this are in contradiction with ours suppositions.\\
{\bf Step 2.} Let $x\in(x_N,\infty)$. Then
$$
h(x)=h(c_N)\frac{b(c_N)}{b(x)}e^{-\int\limits_{c_N}^{x}\!\frac{1}{b(s)}\:ds}\quad,
$$
where we deduce that
$$
\int\limits_{x_N}^{+\infty}\!h(x)\:dx=h(c_N)b(c_N)\int\limits_{x_N}^{+\infty}\frac{1}{b(x)}
e^{-\int\limits_{c_N}^{x}\!\frac{1}{b(s)}\:ds}\:dx.
$$
Let
$$
u(x)=\int\limits_{c_N}^x\!\frac{1}{b(s)}\:ds.
$$
Because $h(c_N)\neq 0$, we have
\begin{itemize}
\item if $b(x)<0$ over $(x_N,+\infty)$, then $u(+\infty)=-\infty$ and 
$$
\int\limits_{x_N}^{\infty}\!h(x)\:dx=h(c_N)b(c_N)\int\limits_{u(x_N)}^{u(\infty)}\!e^{-u}\:du=\infty\mbox{ or }-\infty
$$
\item if $b(x)>0$ for all $x\in(x_N,\infty)$, then
$$
\lim_{x\searrow x_N}b(x)h(x)=h(c_N)b(c_N)e^{\int\limits_{x_N}^{c_N}\!\frac{1}{b(s)}\:ds}=\infty\mbox{ or }-\infty
$$
\end{itemize}
All this are again in contradiction with ours suppositions.\\
{\bf Step 3.} The case where $x\in (-\infty,x_1)$ can be trated like the Step 2.\\
{\it Neccesity.} Suppose in contrary that one of
$$
\int\limits_{-\infty}^{c_0}\!\frac{1}{b^{+}(x)}\:dx\quad\mbox{ or }\quad\int\limits_{c_N}^{+\infty}\!\frac{1}{b^{-}(x)}\:dx
$$
is finite. We work only in the case where
$$
\int\limits_{c_N}^{+\infty}\!\frac{1}{b^{-}(x)}\:dx<\infty
$$
and the other case can be treated in the same way. Define
$$
h(x)=\left\lbrace\begin{array}{lcl}
\frac{b(c_N)}{b(x)}e^{-\int\limits_{c_N}^{x}\!\frac{1}{b(s)}\:ds}&,&x>x_N\\
0&,&x\leq x_N
\end{array}
\right.
$$
We have $-(bh)^{'}=h$ on $(-\infty,x_N)$ and $(x_N,\infty)$. Since
$$
\lim_{x\searrow x_N}b(x)h(x)=b(c_N)e^{\int\limits_{x_N}^{c_N}\!\frac{1}{b(s)}\:ds}=0, 
$$
the function $h$ is again a solution of $-(bh)^{'}=h$ in the sense of distribution, which is in contradiction with the fact that $\cal A$ is $\left(L^\infty\left(\R,dx\right)\right)$-unique. \fin\\
In the multidimensional case $d\geq 2$, the main result of this section is
\begin{thm}\label{3.6}
Let $b:\R^d\rightarrow\R^d$ be a function of the class $C^1(\R^d)$ such that $b(x)\neq 0$ for all $|x|\geq R$. Suppose that there is a locally bounded function $\beta:\R^{+}\rightarrow\R$ such that
$$
\left(b(x)\frac{x}{|x|}\right)^{-}\leq\beta(|x|)\quad\forall|x|\geq R\quad.
$$
If
$$
\int\limits_R^\infty\!\frac{1}{\beta(x)}dx=\infty\quad,
$$
then $\left({\cal A},C_0^\infty(\R^d)\right)$ is
$\left(L^\infty\left(\R^d,dx\right),{\cal
C}\left(L^\infty,L^1\right)\right)$-unique. In particular, for all $h\in L^1(\R^d,dx)$, the mass transport equation
$$
\left\lbrace\begin{array}{l}
\partial_t\rho(t,x)=-div(b\rho(t,x))\\
\rho(0,x)=h(x)
\end{array}
\right.
$$
has one $L^1(\R^d,dx)$-unique weak solution.
\end{thm}
\dem For all $x\in\R^d$, consider $(X_t(x))_{0\leq t<e}$, where $e$ is the explosion time, the solution of the equation
$$
\left\lbrace\begin{array}{l}
\frac{dX_t}{dt}=b(X_t)\\
X(0)=x
\end{array}
\right..
$$
Then the family $\{P_t\}_{t\geq 0}$,
$$
P_t(x)=f(X_t(x))
$$
is a $C_0$-semigroup on $L^\infty(\R^d,dx)$ with respect to the topology ${\cal C}(L^\infty,L^1)$.\\
{\bf Step 1.} We first prove that if $f\in C_0^1$, the there exists $(f_n)_{n\in\N}\subset C_0^\infty(\R^d)$ such that $f_n\rightarrow f$ and ${\cal A}f_n\rightarrow{\cal A}f$ in the topology ${\cal C}(L^\infty,L^1)$.\\
Indeed, let $supp f\subset B(0,N)$. Then by convolution method, there exists $(f_n)_{n\in\N}\in{\cal C}_0^\infty(\R^d)$ such that $supp\:f_n\subset B(0,N+1)$, for all $n\geq 1$, $f_n\rightarrow f$ and $\nabla f_n\rightarrow\nabla f$ uniformely over $\R^d$. Thus
$$
b\nabla f_n\rightarrow b\nabla f
$$
uniformely over $\R^d$.\\
{\bf Step 2.} It remains to prove that $C_0^1$ is a core for $\cal A$. To that purpose, by \cite[Lemma 2.4, p.572]{wu-zhang'06}, it is enough to show that 
$$
P_t{\cal C}^1_0\subset{\cal C}^1_0
$$
or, equivalently, to establish 
$$
\lim_{|x|\rightarrow\infty}|X_t(x)|=\infty.
$$
Consider
$$
\tau_n=\inf\{t\:|\:|X_t|=n\}
$$
and
$$
\tau_\infty=\lim_{n\rightarrow\infty}\tau_n=e.
$$
For all $t<\tau_R\wedge\tau_\infty$, we have
$$
\frac{d|X_t(x)|}{dt}=\frac{X_t(x)}{|X_t(x)|}b(X_t(x))\geq-\beta(|X_t(x)|).
$$
Let
$$
h(x)=\int\limits_{R}^{x}\!\frac{1}{\beta(s)}\:ds.
$$
Then we have
$$
\frac{d}{dt}h(|X_t(x)|)=\frac{1}{\beta(|X_t(x)|)}\frac{X_t(x)}{|X_t(x)|}b(X_t(x))\geq -1
$$
where it follows that
$$
h(|X_t(x)|)\geq h(|x|)-t\quad,\quad\forall t\in[0,\tau_R\wedge\tau_\infty).
$$
Consequently
$$
\lim_{|x|\rightarrow\infty}|X_t(x)|=\infty
$$
where we deduce that $\left({\cal A},C_0^\infty(\R^d)\right)$ is
$\left(L^\infty\left(\R^d,dx\right),{\cal
C}\left(L^\infty,L^1\right)\right)$-unique. \fin

\bibliographystyle{plain}

\end{document}